\documentclass[12pt]{article}

\usepackage{amsthm,amssymb,latexsym,graphics,amscd,amsmath,amscd,enumerate,color}

\usepackage[english]{babel}

\newcommand{\Cl}{\mathcal{C}\ell}

\newcommand{\mb}{\mathbb}

\newcommand{\Sc}{\mb{S}_{\mb{C}}}
\newcommand{\Sp}{\mbox{Spin}}
\newcommand{\Spc}{\mbox{Spin}^{\mb{C}}}

\newcommand{\vol}{\mbox{vol}}
\newcommand{\rk}{\mbox{rk}}

\newtheorem{prop}{Proposition}

\newtheorem{teo}{Theorem}
\newtheorem*{mthm}{Main Theorem}

\begin{document}

\title{(1,1)-forms acting on Spinors on K\"ahler Surfaces}
\author{Rafael de Freitas Le\~ao}
\maketitle

\begin{abstract}
It is known that, for Dirac operators on Riemann surfaces twisted by line bundles with Hermitian-Einstein connections, it is possible to obtain estimates for the first eigenvalue in terms of the topology of the twisting bundle \cite{JL2}. Attempts to generalize topological estimates for higher rank bundles or higher dimensional manifolds have been so far unsuccessful. In this work we construct a class of examples which indicates one problem that arises on such attempts to derive topological estimates.
\end{abstract}

\section{Introduction}

Let $(M,g,J)$ be a K\"ahler manifold of complex dimension $n$, and let $E \rightarrow M$ be a holomorphic Hermitian vector bundle over $M$ with connection $\nabla^A$ compatible with the holomorphic and hermitian structures (the Chern connection). Using the complex structure of $M$, this connection can be decomposed as $\nabla^A = \partial_A + \bar{\partial_A}$, and we can consider the associated twisted Dolbeault Laplacian:
\begin{displaymath}
  \Delta_{\bar{\partial}}=\bar{\partial}_A \bar{\partial}_A^* +  
  \bar{\partial}_A^* \bar{\partial}_A.
\end{displaymath}

When restricted to sections of $E$, the Dolbeault Laplacian simplifies to $\Delta_{\bar{\partial}}=\bar{\partial}_A^* \bar{\partial}_A$ and the K\"ahler identities for the connection $\nabla^A$ \cite{DK} can be used to relate the Dolbeault Laplacian to the connection Laplacian ${\nabla^A}^* \nabla^A$:
\begin{equation} \label{eq_kah_id}
  \Delta_{\bar{\partial}}\mid_{\Omega^{(0,0)(E)}} = 
    \bar{\partial}_A^* \bar{\partial}_A = \frac{1}{2} {\nabla^A}^* 
    \nabla^A - \frac{i}{2} \Lambda F_A,
\end{equation}
where $\Lambda F_A$ is the contraction of the curvature 2-form $F_A$ by the K\"ahler form $\omega$.

In some cases, the term $\Lambda F_A$ simplifies  and the above equation leads to estimates for the eigenvalues of the Dolbeault Laplacian. One case where this happen is when the connection $\nabla^A$ is Hermitian-Einstein\footnote{Recall that a  connection $\nabla^A$ is Hermitian-Einstein if $\Lambda F_A = c \mb{I}_E$, where $c$ is a topological constant given by
\begin{displaymath}
  \frac{2 \pi \deg(E)}{(n-1)! \rk(E) \vol(M)}.
\end{displaymath}}. In this case $\Lambda F_A$ is proportional to the identity and the K\"ahler identity (\ref{eq_kah_id}) can be used to obtain the following lower bound for the eigenvalues of the Dolbeault Laplacian on sections of $E$:
\begin{displaymath}
  \lambda \geq \frac{- \pi \deg(E)}{(n-1)! \rk(E) \vol(M)}.
\end{displaymath}

It was shown in \cite{JL2} that (\ref{eq_kah_id}) can also be derived as a particular case of a convenient Weitzenb\"ock formula. With this, it was possible to use twistor techniques to improve the initial estimate thus obtaining
\begin{equation} \label{eq_sharp}
  \lambda \geq -\frac{2n}{2n-1} \frac{\pi \deg(E)}{(n-1)! \rk(E) 
    \vol(M)}.
\end{equation}

In the same article, we got the proper Weitzenb\"ock formula using the identification between spinors and differential forms on a K\"ahler manifold along with the identity $D_A^2 = 2 \Delta_{\bar{\partial}}$ obtained from this identification. Then the Weitzenb\"ock formula for $D_A$ recovers the above K\"ahler identity. On Riemann surfaces this also provides an estimate for $D_A$ and the explicit computations of \cite{AP} shows that (\ref{eq_sharp}) is sharp.

It is iteresting to note that the relation between the Dirac operator and the Dolbeault Laplacian was used only for sections of a bundle $E$ over a Riemann surface, although it is valid in general. This leads to the natural question whether it is possible to explore this relation in more general cases.

The initial attempt to use the relation between the Dolbeault Laplacian and the Dirac operator fails in higher dimensions because on this situation we need to know estimates for the eigenvalues of the Dolbeaul Laplacian restricted to (0,p)-forms with values on $E$, $\Omega^{0,p}(E)$.

However, on sections of $\Omega^{(0,p)}$ the Weitzenb\"ock formula becomes more complicated. The complication appears because the term involving the curvature $F_A$ cannot, in general, be directly related to the topology of $E$. In the present article we show the following result, which shows one possible reason for this:

\begin{mthm}
Let $E \rightarrow M$ be a holomorphic line bundle over a K\"ahler surface\footnote{By a K\"ahler surface we understand a K\"ahler manifold of complex dimension 2.} $(M,g,J)$. Let $\nabla^A$ be any connection on $E$ compatible with the holomorphic structure and $F_A$ the curvature 2-form of $\nabla^A$. Then, as an operator acting on spinors, $F_A$ is indefinite for every $p \in M$ such that $F_A(p) \neq 0$.
\end{mthm}

This shows that, in general, an estimate along the lines of \cite{Ba} is the best possible. Furthermore, this also shows that attempts to obtain estimates for the eigenvalues of the twisted Dirac operator in higher rank bundles must investigate if $F_A$ can be made into a definite operator.

The proof of this result will be carried in two cases. First we consider anti-selfdual connections and, after that, the more straightforward case of selfdual connections.

\section{Anti-Selfdual $U(1)$ connections}

Let $M$ be a K\"ahler manifold with complex dimension $2$. All complex manifolds carry a canonical $\Spc$-structure, and in this structure the spinor bundle is explicitly described in terms of forms:
\begin{displaymath}
  \begin{split}
    \mb{S}_{\mb{C}} &\simeq \wedge^{0,*} M = \oplus_{i=0}^2 
    \wedge^{0,i} M, \\
    \mb{S}_{\mb{C}}^+ &\simeq \oplus_{i \hspace{1ex} 
      even}\wedge^{0,i} M, \\
    \mb{S}_{\mb{C}}^- &\simeq \oplus_{i \hspace{1ex} odd} 
      \wedge^{0,i} M.
  \end{split}
\end{displaymath}
Consequently, the twisted case is described by
\begin{displaymath}
  \mb{S}_{\mb{C}} \otimes E \simeq \wedge^{0,*} M \otimes E = 
    \Omega^{0,*}(E).
\end{displaymath}

This description is very useful, mainly because of two reasons. First, we can explicitly describe the action of $\Cl(T^*M)$ on $\Sc$. For this, consider an adapted frame $\{ \xi^i, \bar{\xi}^i \}$ of $T^*M \otimes \mb{C}$. Then, in this frame, the Clifford action is given by:
\begin{equation} \label{cliff_act}
  \begin{split}
    \xi^i \cdot &= - \sqrt{2} \bar{\xi}^i \lrcorner, \\
    \bar{\xi}^i \cdot &= \sqrt{2} \bar{\xi}^i \wedge.
  \end{split}
\end{equation}
Secondly, the twisted Dirac operator can be described in terms of Cauchy-Riemann operators: if $\nabla^A$ is a connection on $E \rightarrow M$, the complex structure of $M$ produces the splitting
\begin{displaymath}
  \begin{split}
    \nabla^A = \partial_A + \bar{\partial}_A, \\
    \partial_A: \Omega^{p,q}(E) \rightarrow \Omega^{p+1,q}(E), \\
    \bar{\partial}_A:\Omega^{p,q}(E) \rightarrow \Omega^{p,q+1}(E),
  \end{split}
\end{displaymath}
and the twisted Dirac operator is given by
\begin{displaymath}
  D_A = \sqrt{2} \left( \partial_A + \bar{\partial}_A \right).
\end{displaymath}
The case of the twisted Dirac operator associated with a $\Sp$-structure can also be described by these identifications; we only must remember that the two spinor spaces are related by $\mb{S}_{\mb{C}} = \mb{S} \otimes k^{-\frac{1}{2}}$, where $k = \wedge^{0,n} M$.

Another important fact about K\"ahler manifolds with complex dimension 2 is that the 2-forms decompose in self-dual forms, $\Omega^+$, and anti-self-dual forms, $\Omega^-$, and that
\begin{equation} \label{eq_dec}
  \begin{split}
    \Omega^+ &= \Omega^{2,0} \oplus \Omega^0 \omega \oplus 
      \Omega^{0,2}, \\
    \Omega^- &= \Omega^{1,1}_0,
  \end{split}
\end{equation}
where $\omega$ is the K\"ahler form and $\Omega^{1,1}_0$ is the space of (1,1)-forms orthogonal to $\omega$ \cite{DK}.

Using the adapted frame $\{ \xi^i, \bar{\xi}^i \}$ we can explicitly describe the action of elements of $\Omega^{\pm}$ on spinors. First, note that the K\"ahler form can be written as
\begin{displaymath}
  \omega = i \left( \xi^1 \wedge \bar{\xi}^1 + \xi^2 \wedge 
    \bar{\xi}^2 \right),
\end{displaymath}
and a basis for $\Omega^{1,1}_0$  is given by $\{ \xi^1 \wedge \bar{\xi}^2,\xi^2 \wedge \bar{\xi}^1,\xi^1 \wedge \bar{\xi}^1-\xi^2 \wedge \bar{\xi}^2 \}$. Therefore, if $F_A \in \Omega^-$, locally we can write
\begin{displaymath}
  F_A = a \xi^1 \wedge \bar{\xi}^2 + b\xi^2 \wedge \bar{\xi}^1 + 
    c\left( \xi^1 \wedge \bar{\xi}^1-\xi^2 \wedge \bar{\xi}^2 
    \right).
\end{displaymath}

\begin{prop} \label{prop_act}
If $F_A \in \Omega^-$, the action of $F_A$ on $\mb{S}^-$ is given by
\begin{displaymath}
  F_A = 2
  \begin{pmatrix}
    c & b \\
    a & -c
  \end{pmatrix}.
\end{displaymath}
\end{prop}

{\bf Proof.} A 2-form $\alpha \wedge \beta$ acts on spinors through Clifford multiplication by means of the identification
\begin{displaymath}
  \alpha \wedge \beta \simeq \frac{1}{2}\left( \alpha \beta - \beta 
    \alpha \right).
\end{displaymath}
Using the action described in (\ref{cliff_act}) we calculate
\begin{displaymath}
  \begin{split}
    \xi^1 \wedge \bar{\xi}^2 \cdot \psi &= \frac{1}{2} \left( \xi^1
      \bar{\xi}^2 - \bar{\xi}^2 \xi^1 \right) \cdot \psi \\
    &= \frac{1}{2} \left[ \xi^1 \cdot \left( \bar{\xi}^2 \cdot 
      \psi \right) - \bar{\xi}^2 \cdot \left( \xi^1 \cdot \psi 
      \right) \right] \\
    &= - \bar{\xi}^1 \lrcorner \left( \bar{\xi}^2 \wedge \psi 
      \right) + \bar{\xi}^2 \wedge \left( \bar{\xi}^1 \lrcorner 
      \psi \right).
  \end{split}
\end{displaymath}
For 4-manifolds $\mb{S}_{\mb{C}}^-$ is just $\Omega^{0,1}(E)$; so if $\psi \in \mb{S}_{\mb{C}}^-$ we have\footnote{Strictly, elements of $\Omega^{(0,1)}(E) \simeq \Gamma(E) \otimes \wedge^{(0,1)} M$ are of the form $\psi = \psi_1 \otimes \bar{\xi}^1 + \psi_2 \otimes \bar{\xi}^2$ and the Clifford action is given by $c(\alpha) (\psi_i \otimes \bar{\xi}^i) = \psi_i \otimes ( c(\alpha) \bar{\xi}^i)$. So, to simplify notation, we just write $\psi_i \otimes \bar{\xi}^i \sim \psi_i \bar{\xi}^i$, and the action is as written.}
\begin{displaymath}
  \psi = \psi_1 \bar{\xi}^1 + \psi_2 \bar{\xi}^2,
\end{displaymath}
and the above expression simplifies to
\begin{displaymath}
  (\xi^1 \wedge \bar{\xi}^2) \cdot \psi = \psi_1 \bar{\xi}^2 +
    \psi_1 \bar{\xi}^2 = 2 \psi_1 \bar{\xi}^2
\end{displaymath}
or, in matrix form,
\begin{displaymath}
  (\xi^1 \wedge \bar{\xi}^2) \cdot \psi =
  \begin{pmatrix}
    0 & 0 \\
    2 & 0
  \end{pmatrix} \cdot
  \left[
  \begin{array}{c}
    \psi_1 \\
    \psi_2
  \end{array}
  \right].
\end{displaymath}
The other terms are calculated in the same manner and are given, in matrix form, by
\begin{displaymath}
  \begin{split}
    (\xi^2 \wedge \bar{\xi}^1) \cdot \psi =
    \begin{pmatrix}
      0 & 2 \\
      0 & 0
    \end{pmatrix} \cdot
    \left[
    \begin{array}{c}
      \psi_1 \\
      \psi_2
    \end{array}
    \right], \\
    \left( \xi^1 \wedge \bar{\xi}^1 - \xi^2 \wedge \bar{\xi}^2 
      \right) \cdot \psi =
    \begin{pmatrix}
      2 & 0 \\
      0 & -2
    \end{pmatrix} \cdot
    \left[
    \begin{array}{c}
      \psi_1 \\
      \psi_2
    \end{array}
    \right].
  \end{split}
\end{displaymath}
and the result follows. $\square$

To complete the characterization we need to know how $F_A$ acts on $\mb{S}^+$.
\begin{prop} \label{prop:cli_act}
The Clifford action of an $(1,1)$-form $\alpha$ on sections of \break $\wedge^{0,n}(M) \otimes E$, $\psi_n \in \Omega^{0,n}(E)$, is explicit given by
\begin{equation}
  \alpha \cdot \psi_n = i ( \Lambda \alpha ) \psi_n,
\end{equation}
where $\Lambda \alpha = \omega \lrcorner \alpha$ is the contraction of $\alpha$ by the K\"ahler form $\omega$.
\end{prop}

{\bf Proof.} In \cite[Proposition1]{JL2} it was proved that the action of an (1,1)-form, $\alpha$, on sections of $E$ is given by $-i ( \Lambda \alpha )$. The same techniques can be used to explicitly obtain the result. $\square$

With this we can prove:

\begin{teo} \label{teo_asd}
If $F_A \in \Omega^-$ then, on points $p \in M$ such that $F_A(p) \neq 0$ as a 2-form, $F_A(p)$, as an operator on $\mb{S}_p$,  is indefinite.
\end{teo}

{\bf Proof.} By the above proposition, the action of a (1,1)-form $\alpha$ on $\wedge^{0,0} M$ and $\wedge^{0,2} M$ is given by
\begin{displaymath}
  \alpha \cdot \psi = \mp i ( \Lambda \alpha ) \psi,
\end{displaymath}
where the sign is minus on $\wedge^{0,0} M$ and plus on $\wedge^{0,2} M$.

Then, if $F_A \in \Omega^-$, the decomposition (\ref{eq_dec}) implies that $\Lambda F_A = 0$, and $F_A$ acts trivially on $\mb{S}_{\mb{C}}^+$. The only non trivial part is the action of $F_A$ on $\mb{S}_{\mb{C}}^-$, which is given by proposition (\ref{prop_act}).

With the same notation as in the previous proposition, and knowing that the action of $F_A$ on $\mb{S}_{\mb{C}}$ is Hermitian \cite{BL} we find that the eigenvalues of $F_A$, as an operator on $\mb{S}_{\mb{C}}$, are
\begin{displaymath}
  \{ 0, \sqrt{c^2+ab}, -\sqrt{c^2+ab} \}.
\end{displaymath}

Because the representation of $\Cl(TM)$ on $\mb{S}_{\mb{C}}$ is faithful we conclude that for points $p \in M$ where $F_A(p) \neq 0$ the eigenvalues $\pm \sqrt{c^2+ab}$ cannot be zero, so $F_A$ is indefinite. $\square$

\section{Selfdual $U(1)$ connections}

For compatible self-dual connections the decomposition (\ref{eq_dec}) implies that the curvature has the form $F_A = f \omega$, where $\omega$ denotes the K\"ahler form and $f$ is a function on $M$. 

Using proposition (\ref{prop:cli_act}) and \cite[Proposition1]{JL2} we have:
\begin{prop} \label{prop_sd}
If $F_A$ is of type (1,1), the action of $F_A$ on $\mb{S}^+ \simeq \wedge^{0,0} M \oplus \wedge^{0,2} M$ is given by
\begin{displaymath}
  F_A = i
  \begin{pmatrix}
    - \Lambda F_A & 0 \\
    0 & \Lambda F_A
  \end{pmatrix}.
\end{displaymath}
\end{prop}

Using this we have:

\begin{teo} \label{teo_sd}
If $F_A \in \Omega^+$ then, on points $p \in M$ such that $F_A(p) \neq 0$ as a 2-form, $F_A(p)$ as an operator on $\mb{S}_p$  is indefinite.
\end{teo}

{\bf Proof.} Using the calculations of proposition (\ref{prop_act}) we can explicitly verify that $\omega$ acts as a null operator on $\mb{S}^-$. Thus $F_A$ acts trivially on $\mb{S}^-$ and the action of $F_A$ on $\mb{S}^+$ is given by the above proposition. Therefore, for $F_A \neq 0$, is indefinite. $\square$

Combining theorems (\ref{teo_asd}) and (\ref{teo_sd}), and decomposition (\ref{eq_dec}), we obtain the main theorem.

\end{document}